\documentclass[
11pt,
authoryear,
times
]{elsarticle}

\usepackage{graphicx,color}

\usepackage[colorlinks=true, breaklinks=true, pdfstartview=FitV, linkcolor=red, citecolor=blue, urlcolor=black]
{hyperref}

\usepackage[T1]{fontenc}

\usepackage{algorithm}
\usepackage{algorithmicx}
\usepackage{algpseudocode}

\newcommand{\Vec}[1]{\mbox{\boldmath$#1$}}
\journal{Celestial Mechanics and Dynamical Astronomy}

\begin{document}

\begin{frontmatter}

\title{A Hopf variables view on the libration points dynamics.}

\author[roa]
{ Martin Lara\fnref{SDG} 
}
\ead{mlara0@gmail.com}


\fntext[SDG]{ GRUCACI, University of La Rioja, and Space Dynamics Group -- UPM  }

\address[roa] {C/ Luis de Ulloa, s.n., 26004 Logro\~no, Spain}

\date{}
\begin{abstract}
The dynamics about the libration points of the Hill problem is investigated analytically. In particular, the use of Lissajous variables and perturbation theory allows to reduce the problem to a one degree of freedom Hamiltonian depending on two physical parameters. The invariant manifolds structure of the Hill problem is then disclosed, yet accurate computations are limited to energy values close to that of the libration points.
\end{abstract}

\begin{keyword}
libration points \sep Hill's problem \sep center manifold \sep perturbation theory \sep elliptic oscillator \sep periodic orbits \sep Lissajous variables \sep Hopf variables \sep resonant normal form


\end{keyword}

\end{frontmatter}

\section{Introduction}

The Hill problem is an useful approximation of the restricted three-body problem ---not to be confused with a particular case \citep{HenonPetit1986}--- that, further than its original application to the computation of the orbit of the moon around the earth \citep{Hill1878}, can be representative of the dynamics of an object under the gravitational atraction of different pairs of solar system bodies. Indeed, when using suitable units of length and time the Hill problem does not depend on any parameter, and, therefore, its application to different scenarios becomes a simple matter of scaling \citep{Szebehely1967}. In particular, the Hill problem is well suited by itself to study the dynamics about asteroids, yet it may need to be amended to include the important effect of the solar radiation pressure \citep[see][and references therein]{GarciaScheeresMacIness2015}. But it can be used too in the description of the most relevant features of the dynamics around planetary satellites \citep{LaraRusselVillacl2007Meccanica,LaraPalacianRussell2010}, a case that may require to further superimpose to the third-body dynamics the nonspherical disturbances of the central body, which can notably modify the orbital behavior close to the origin \citep[see][for instance]{LidovYarskaya1974,Vashkovyak1996,ScheeresGumanVillac2001,LaraSanJuan2005,RussellLara2009}. The Hill problem equations are useful also in the investigation of satellite encounters \citep{PetitHenon1986}, and can capture the bulk of the dynamics of coorbital motion, with different applications to relative spacecraft motion \cite[see][and references therein]{KasdinGurfilKolemen2005}.
\par

The Hill problem has been thoroughly studied numerically by the propagation of periodic and quasi-periodic orbits \citep{Henon1969,Henon1970,Henon1974,Henon2003,Michalodimitrakis1980}, as well as escape trajectories \citep{VillacScheeres2003}. Other studies provide detailed accounts of the dynamics, including the global description of the planar case for values of the energy corresponding to bounded motion \citep{SimoStuchi2000}. Besides, due to its interest in spacecraft mission design, special emphasis has been given to the study of the stable and unstable manifolds associated to Lissajous orbits, which can be effectively computed from the investigation of the center manifold of the collinear libration points \citep{GomezMarcoteMondelo2005,Masdemont2005}.
\par

The global dynamics of the Hill problem must be necessarily investigated numerically, although purely analytical approaches may provide useful information in those regions of phase space in which the motion can be considered a perturbation of the Keplerian motion ---the normalized solution been commonly constrained to the close vicinity of the primary \citep{SanJuanLaraFerrer2006,Lara2008,LaraPalacianYanguasCorral2010}. On the other hand, the normal form approach is not restricted to the case of perturbed Keplerian motion and is customarily used in the computation of the center manifold of the libration points \citep{GomezJorbaMasdemontSimo1991}. This two degrees of freedom manifold is investigated numerically with the usual tools of non-linear dynamics, as Poincar\'e surfaces of section and the continuation of periodic and quasi periodic orbits \citep{GomezMarcoteMondelo2005}. Optionally, analytical approximations to the existing periodic orbits can be obtained  with the Lindstedt-Poincar\'e method \citep{ZagourasMarkellos1985} \cite[see, also][]{FarquharKamel1973,Richardson1980b}.
\par

Alternatively to the use of Poincar\'e surfaces of section, the dynamics of the center manifold can be approached analytically, at least for energy values close enough to the energy of the libration points. Indeed, the center manifold Hamiltonian of the Hill problem has the form of a two degrees of freedom perturbed harmonic oscillator in the quasi-resonance condition, which is easily cast into the form of a perturbed elliptic oscillator by the standard introduction of a detuning parameter \citep{Henrard1970}. Then, the Hamiltonian is rearranged in the form of an unperturbed term in the 1-1 resonance condition, whereas the terms that have been decoupled with the detuning parameter are incorporated into the perturbation. The dynamics of these classes of resonant systems can be efficiently approached analytically using perturbation theory \citep[see the recent review in][]{MarchesielloPucacco2016} and has been lately applied to the computation of analytical approximations of the Lissajous and Halo orbits in the restricted three-body problem \citep{CellettiPucaccoStella2015}. The former was approached by a standard double normalization of the center manifold Hamiltonian in harmonic-type variables, whereas the computation of the later required a preprocessing of the center manifold Hamiltonian in order to apply resonant perturbation theory \citep{FerrazMello2007}.
\par

On the other hand, the Lissajous transformation \citep{Deprit1991} comes out as a convenient option to the customary use of harmonic variables in dealing with elliptic oscillators. It  was specifically devised to deal with perturbed elliptic oscillators, and reveals particularly well suited to the construction of a resonant normal form Hamiltonian by standard averaging over the elliptic anomaly. The normalized Hamiltonian is of one degree of freedom and, after reformulation in the \citet{Hopf1931} coordinates, provides a complete description of the reduced dynamics on the sphere \citep{DepritElipe1991,Miller1991}.
\par

In the present research it is shown that application of the Lissajous transformation to the center manifold Hamiltonian of the Hill problem allows for the straightforward construction of an analytical solution for the motion about the libration points. The theory only accommodates the cubic and quartic terms of the perturbation expansion of the Hill problem Hamiltonian about the libration points, leading to an extremely simple normalized Hamiltonian. The standard transformation that performs the reduction of the center manifold is provided up to quadratic corrections, and only short-period terms of comparable accuracy are computed in the normalization of the center manifold Hamiltonian in Lissajous variables. Therefore, the practical application of the current theory is constrained to the lower orders of the energy for which these early truncations made sense. Even so, the insights provided by this simple approach go much further than expected and the solution is able to capture the main features of the dynamics about the libration points. Indeed, it not only shows the existence of the planar and vertical Lyapunov orbits, which exist for all values of the energy above the energy of the libration points; but it also shows the main bifurcations of these fundamental orbits, which occur for energy values considerably far away from that of the libration points. Namely, the bifurcation which gives rise to Halo orbits, and the bifurcation and termination of the two-lane bridge of periodic orbits that connects the families of planar and vertical Lyapunov orbits. The construction of a higher order theory, which will notably improve the accuracy of the solution, is just a matter of mechanizing computations and is not discussed here. 
\par

The paper is organized as follows. First of all, basic facts of the Hill problem, including information about the main families of periodic orbits related to the libration point dynamics, are recalled in Section \ref{se:basics}. Next, the construction of the perturbation solution for the motion in the vicinity of the libration points, which consists of the reduction to the center manifold and the consequent removal of short period effects, is approached is Section \ref{se:peso}. It follows the discussion of the reduced phase space in Section \ref{se:reduflow}, where the equilibria of the reduced dynamics are identified with the main existing families of periodic orbits about the libration points. Finally, some validation tests of the proposed solution are presented in  Section \ref{se:validation}.

\section{Hill problem dynamics} \label{se:basics}

In a rotating frame with velocity $\Vec{N}$ in the $z$ axis direction, the $x$ axis defined by the line joining the primaries, the $y$ axis completing a direct frame, and taking one of the primaries as the origin, the Hill problem is defined by the Hamiltonian 
\begin{equation} \label{HillHam}
\mathcal{J}=\mbox{$\frac{1}{2}$}(\Vec{P}\cdot\Vec{P})-\Vec{N}\cdot(\Vec{p}\times\Vec{P})-\Omega(\Vec{p}),
\end{equation}
where $\Vec{p}\equiv(p_x,p_y,p_z)$ is position, its conjugate momentum $\Vec{P}\equiv(P_x,P_y,P_z)$ is velocity in the inertial frame, and $\Omega$ is the potential function
\[
\Omega=\frac{\mu}{R}-\frac{N^2}{2}(R^2-3p_x^2),
\]
where $R=\|\Vec{p}\|$ and $\mu$ is the gravitational parameter. 
\par

After scaling units of length by $\mu^{1/3}$ and time by $1/N$, the Hill problem Hamiltonian is rewritten
\begin{equation} \label{Jota}
\mathcal{J}=\frac{1}{2}\left(P_x^2+P_y^2+P_z^2\right)+P_x p_y-p_x P_y-\frac{1}{R}+\frac{1}{2} \left(R^2-3 p_x^2\right),
\end{equation}
showing that the Hill problem does not depend on any parameter. In these non-dimensional units, the Hamiltonian equations stemming from Eq.~(\ref{Jota}) are
\begin{eqnarray} \label{dpx}
\dot{p}_{x} &=& P_x+p_{y}, \\
\dot{p}_{y} &=& P_y-p_{x}, \\
\dot{p}_{z} &=& P_z, \\
\dot{P}_{x} &=& -\frac{1}{R^3}p_x+2p_x+P_y, \\
\dot{P}_{y} &=& -\frac{1}{R^{3}}p_y-p_y-P_x, \\ \label{dPz}
\dot{P}_{z} &=& -\frac{1}{R^{3}}p_z-p_z,
\end{eqnarray}
where over dots denote derivatives in the rotating frame. From Eqs.~(\ref{dpx})--(\ref{dPz}) it is immediately apparent that planar motions $p_z=P_z=0$ exist, as well as the two equilibria $\mathcal{L}_{1,2}=\pm(\rho,0,0,0,\rho,0)$, where
\begin{equation} \label{HillRadius}
\rho=3^{-1/3},
\end{equation}
is called the Hill sphere radius, or Hill radius in short. The two equilibria $\mathcal{L}_{1,2}$ are customarily named libration points. Due to the symmetries of the Hill problem with respect to the plane $x=0$, it is enough to discuss the dynamics about just one of the libration points, say $\mathcal{L}_1$. 
\par

The study of the linearized dynamics about the libration points shows that even though they are unstable equilibria, periodic motion originates from them in the form of small vertical oscillations through them, and planar oscillations around them \citep{Szebehely1967}. The so called vertical and planar Lyapunov orbits are then grouped into natural families of periodic orbits which are parameterized by the energy, and are customarily computed by numerical continuation techniques \cite[see][for instance]{Doedeletal2003}. 
\par

\begin{figure}[htbp]
\centering \includegraphics[scale=1, angle=0]{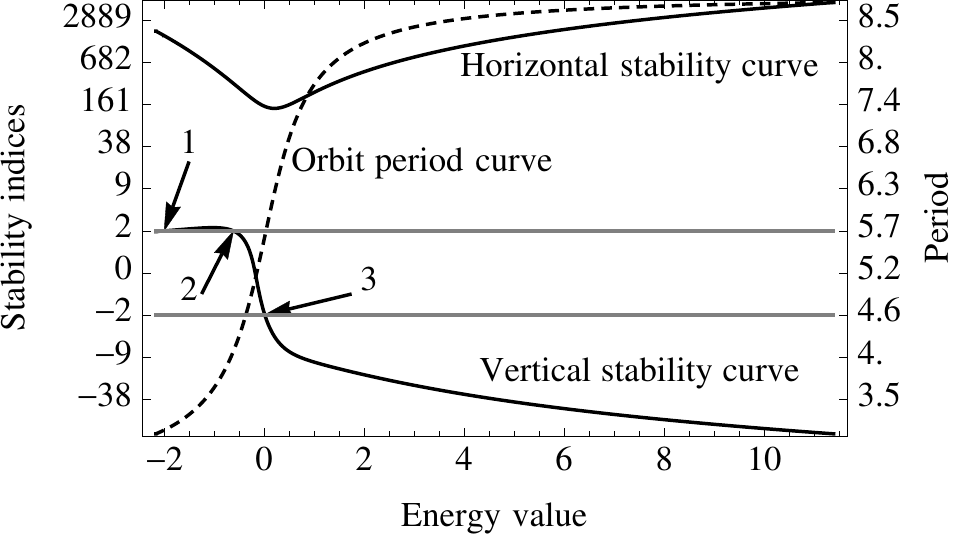}
\caption{Period and stability curves of the family of planar Lyapunov orbits. Arrows point to the vertical bifurcations, which happen at the energy values $\approx-2,-0.6$ and $0$.}
\label{f:planarL}
\end{figure}

The stability of each periodic orbit is characterized by two parameters, say $s_1$ and $s_2$, where orbit stability requires that both indices are real numbers with absolute value less than 2, and bifurcations of new families of periodic orbits may happen when any of the indices crosses this level. The graphic representation of the evolution of these indices along the family provides useful information. Thus, the stability curves of the family of planar Lyapunov orbits is depicted in Fig.~\ref{f:planarL}, where the usual scaling $2\,\mathrm{arcsinh}\,{s}_i/\mathrm{arcsinh}\,2$ is used rather than $s_i$. Because the orbits are planar, one of the indices is related to in-plane perturbations whereas the other is related to out-of-pane perturbations. As shown in the figure, planar Lyapunov orbits are highly unstable, and three crossings of the critical value $|2|$ are observed. The first one occurs close to the energy value $\mathcal{H}(\Vec{p},\Vec{P})\approx-2$, where the Halo orbits emerge; the second crossing occurs for $\mathcal{H}(\Vec{p},\Vec{P})\approx-0.6$, the beginning of a two-lane bridge of periodic orbits which connect the planar and vertical families of Lyapunov orbits. Finally, when $\mathcal{H}(\Vec{p},\Vec{P})\approx0$ a new family of periodic orbits bifurcates with duplication of the period. More details about periodic orbits and other invariant objects of the Hill problem can be consulted in \citep{GomezMarcoteMondelo2005}.
\par

\section{Perturbation solution} \label{se:peso}

First of all, the origin is translated to the libration point by means of the canonical transformation
\begin{equation} \label{T1}
\mathcal{T}_1:(p_x,p_y,p_z,P_x,P_y,P_z)\longrightarrow(x,y,z,X,Y,Z),
\end{equation}
given by
\begin{equation} \label{translation}
p_x=\rho+x, \quad p_y=y, \quad p_z=z, \quad P_x=X \quad P_y=Y+\rho, \quad P_z=Z.
\end{equation}
Then, Eq.~(\ref{HillHam}) is rewritten
\begin{eqnarray} \label{H1ham}
\mathcal{T}_1\circ\mathcal{J}=\mathcal{H} &\equiv& \mbox{$\frac{1}{2}$}(X^2+Y^2+Z^2)-xY+yX+\mbox{$\frac{1}{2}$}(y^2+z^2)-x^2 \\ \nonumber
&& -\frac{x}{\rho^2}-\frac{1}{\sqrt{\left(x+\rho\right)^2+y^2+z^2}}.
\end{eqnarray}
\par

For small values of the distance in Hill units $r=\sqrt{x^2+y^2+z^2}$ compared to the Hill radius $\rho\approx0.7$, the last summand in Eq.~(\ref{H1ham}) can be replaced by the usual expansion in Legendre polynomials. Then, Eq.~(\ref{H1ham}) is written in the form of a perturbation Hamiltonian
\begin{equation} \label{Hillps}
\mathcal{H} =\sum_{n\ge0}\frac{\epsilon^n}{n!}H_n,
\end{equation}
in which the zeroth order term
\begin{equation} \label{H0}
H_0=\mbox{$\frac{1}{2}$}(X^2+Y^2)-(xY-Xy)+2(y^2-2x^2)+\mbox{$\frac{1}{2}$}(Z^2+4z^2),
\end{equation}
is integrable, $\epsilon$ is a formal small parameter indicating the strength of each term $n$ of the perturbation, and
\begin{equation} \label{Hn}
H_n=-\frac{n!}{\rho}\frac{r^{n+2}}{\rho^{n+2}}P_{n+2}(x/r), \qquad (n\ge0),
\end{equation}
where $P_n$ is the Legendre polynomial of degree $n$.

\subsection{Linearized dynamics}

The linear dynamics about the libration points is obtained by truncating the perturbation Hamiltonian (\ref{Hillps}) to retain only the term $H_0$. The last summand of Eq.~(\ref{H0}) matches the Hamiltonian of the simple harmonic oscillator with frequency
\begin{equation} \label{nu}
\nu=2,
\end{equation}
and, therefore, shows that, in the linear approximation, the motion in ${z}$ and ${Z}$ decouples from the rest of the flow about the libration point and comprises small oscillations in the $z$ axis direction which, therefore, is a stable direction in all cases. 
More precisely, the projection of the tangent flow in the $({z},{Z})$ plane is made of ellipses, so that, relative to this plane, the equilibria of the Hill problem are of the center type.
\par

On the other hand, from the Hamilton equations of $H_0$,
\begin{equation} \label{M1}
\left(\begin{array}{c} \dot{x} \\ \dot{y} \\ \dot{X} \\ \dot{Y}\end{array}\right) 
=M_1\left(\begin{array}{c} {x} \\ {y} \\ {X} \\ {Y} \end{array}\right), \qquad
M_1= \left(\begin{array}{cccccc} 
 0 & 1 & 1 & 0 \\ 
-1 & 0 & 0 & 1 \\ 
 8 & 0 & 0 & 1 \\ 
 0 &-4 &-1 & 0 
\end{array}\right),
\end{equation}
is a linear differential system with constant coefficients, whose general solution is made of a linear combination of exponentials
\begin{equation} \label{xi}
\xi_i=\sum_{j=1,4}A_{i,j}\exp(\lambda_jt), \qquad i=1,\dots4,
\end{equation}
where $\xi_i$ stands for $x$, $y$, $X$, and $Y$, respectively, $(A_{i,j})$ is a 4 $\times$ 4 matrix of arbitrary coefficients, which may be expressed as functions of the initial conditions, and the characteristic exponents $\lambda_j$ are the eigenvalues of $M_1$. Namely $\lambda_{1,2}=\pm\lambda$, $\lambda_{3,4}=\pm\omega\jmath$, with $\jmath=\sqrt{-1}$ and
\begin{equation} \label{lw}
\lambda=(2\sqrt{7}+1)^{1/2}, \qquad \omega=(2\sqrt{7}-1)^{1/2}.
\end{equation}
\par

Because $\lambda_{1,2}$ are real, the general solution in Eq.~(\ref{xi}) has a hyperbolic component. On the other hand, the characteristic exponents $\lambda_{3,4}$ are pure imaginary, thus giving place to an elliptic or center-type component. Therefore, the equilibrium points of the Hill problem are of the saddle $\times$ center $\times$ center type. Notably, the same dynamical behavior results from the Hamiltonian in separate variables
\begin{equation} \label{separable}
K_0=\lambda{x}_1X_1+\frac{1}{2}(Y_1^2+\omega^2y_1^2)+\frac{1}{2}(Z_1^2+\nu^2z_1^2),
\end{equation}
which is obtained after the canonical transformation
\begin{equation} \label{T2}
\mathcal{T}_2:(x,y,z,X,Y,Z)\longrightarrow(x_1,y_1,z_1,X_1,Y_1,Z_1),
\end{equation}
given by the linear transformation
\begin{equation} \label{x1tox}
\left(\begin{array}{c} x \\ y \\ X \\ Y \end{array}\right) =A\left(\begin{array}{c} x_1 \\ y_1 \\ X_1 \\ Y_1 \end{array}\right),
\end{equation}
in which
\begin{equation} \label{Agerard}
A=\left(
\begin{array}{cccc}
 2 \lambda/\sigma & 0 & -2 \lambda/\sigma & 2/\tau \\ 
 (\lambda^2-9)/\sigma & -(\omega^2+9)/\tau & (\lambda^2-9)/\sigma & 0 \\
 (\lambda^2+9)/\sigma & (9-\omega^2)/\tau & (\lambda^2+9)/\sigma & 0 \\
 \lambda(\lambda^2-7)/\sigma & 0 & \lambda(7-\lambda^2)/\sigma & -(\omega^2+7)/\tau
\end{array}
\right),
\end{equation}
is obtained based on the eigenvector decomposition of $M_1$ in Eq.~(\ref{M1}), cf.~\citep{JorbaMasdemont1999}. Note that Eq.~(\ref{Agerard}) is slightly different from the equivalent matrix in \citep{GomezMarcoteMondelo2005}, where the last column in Eq.~(\ref{Agerard}) appears multiplied by $\omega$. This is just a consequence of the different form of Eq.~(\ref{separable}), which is intentionally chosen in preparation for a following Lissajous transformation, from the zero order Hamiltonian in \citep{GomezMarcoteMondelo2005}.
\par

\subsection{Reduction to the center manifold}

The flow derived from Eq.~(\ref{separable}) admits the three integrals
\[
J_x=x_1X_1, \qquad J_y=Y_1^2+\omega^2y_1^2, \qquad J_z=Z_1^2+\nu^2z_1^2,
\]
which, for given values of each of them, define corresponding invariant manifolds of the linearized motion. When the motion is constrained to the manifold $J_x=0$ then  $\dot{x}_1=\dot{X}_1=0$ and the saddle component is removed. Therefore, the manifold $J_x=0$ is of the center $\times$ center type, and, for this reason, is called the center manifold.
\par

The existence of the integral $J_x$, and, as a consequence, the center manifold, is not limited to the linear dynamics and can be extended to the nonlinear terms of the transformed Hamiltonian
\begin{equation} \label{Kam}
\mathcal{T}_2\circ\mathcal{H}=\mathcal{K}\equiv\sum_{n\ge0}\frac{\epsilon^n}{n!}K_n(x_1,y_1,z_1,X_1,Y_1,Z_1),
\end{equation}
which is obtained after applying the transformation defined by Eqs.~(\ref{x1tox})--(\ref{Agerard}) to all the summands of Eq.~(\ref{Hillps}).
\par

The procedure for extending the integral $J_x$ to the nonlinear terms consist in finding a canonical transformation
\begin{equation} \label{T3}
\mathcal{T}_3:(x_1,y_1,z_1,X_1,Y_1,Z_1)\longrightarrow(x_2,y_2,z_2,X_2,Y_2,Z_2),
\end{equation}
that converts Eq.~(\ref{Kam}) into a normal form such that, for instance \citep[other possibilities may exist, cf.][]{GomezMarcoteMondelo2005}, in the new variables all the monomials
\[
M_k=Q_kx_2^{m_1}X_2^{m_2}y_2^{m_3}Y_2^{m_4}z_2^{m_5}Z_2^{m_6}, \quad k=(m_1,m_2,m_3,m_4,m_5,m_6),
\]
with $m_1\ne{m}_2$ are removed from the Hamiltonian \citep{Giorgillietal1989}. This yields a transformed Hamiltonian
\begin{equation} \label{KJ}
\mathcal{T}_3\circ\mathcal{K}=\sum_{n\ge0}\frac{\epsilon^n}{n!}K_n(y_2,z_2,Y_2,Z_2;J),
\end{equation}
where $J=x_2X_2$ is an integral. Equation (\ref{KJ}) is in the required normal form, and, therefore, has the center manifold $J=0$.
\par

The transformation $\mathcal{T}_3$ is computed by canonical perturbation theory \citep{Deprit1969}. The construction of the center manifold Hamiltonian is simpler and better understood when using complex variables \citep{Kummer1976}, because they make trivial the solution of the homological equation of the perturbation method. However, this additional change of variables is not necessary for the low order of the theory presented here, and the computations have been made directly in the subindex 2 variables.

After neglecting terms of $\mathcal{O}(\epsilon^3)$ and higher, the Hamiltonian of the center manifold is obtained by making $J=0$, in Eq.~(\ref{KJ}), viz.
\begin{equation} \label{centerHam}
\mathcal{C}=\sum_{n=0}^{2}\frac{\epsilon^n}{n!}C_n,
\end{equation}
with the summands
\begin{eqnarray}
C_0 &=& \frac{1}{2}(Y_2^2+\omega^2y_2^2)+\frac{1}{2}(Z_2^2+\nu^2z_2^2), \\
C_1 &=& \frac{1}{56}\rho^2\tau\left[-\frac{27}{2}y_2^2 -3(2\omega^2-5) z_2^2 +\frac{1}{9}(19-4\omega^2)Y_2^2\right]Y_2, \\
C_2 &=& \rho\left[-\frac{81}{1083488}(1322\omega^2+22707)y_2^4
+\frac{27}{270872}(643\omega^2+22588)y_2^2 Y_2^2 \right. \qquad \\ \nonumber &&
-\frac{27}{812}(\omega^2-16) y_2 Y_2 z_2 Z_2
-\frac{27}{1122184}(36962\omega^2-19773)y_2^2 z_2^2 \\ \nonumber
&& +\frac{27}{1624}(5\omega^2+36)y_2^2 Z_2^2
+\frac{1}{2437848}(82144\omega^2-445831)Y_2^4 \\ \nonumber
&& +\frac{9}{561092}(55909\omega^2-137470)Y_2^2z_2^2
+\frac{3}{812}(\omega^2-16)Y_2^2Z_2^2 \\ \nonumber
&& \left. +\frac{27}{1624}(34\omega^2-225)z_2^4 
+\frac{27}{812}(3\omega^2+10)z_2^2 Z_2^2\right].
\end{eqnarray}
\par

Up to the first order, the transformation is given by $y_1=y_2$, $z_1=z_2$, $Y_1=Y_2$, $Z_1=Z_2$, and
\begin{equation} \label{directCart}
x_1=x_2-\rho^2\sigma(\lambda\Delta_1+\Delta_2),\qquad X_1=X_2+\rho^2\sigma(\lambda\Delta_1-\Delta_2),
\end{equation}
with
\begin{eqnarray} \label{deltax}
\Delta_1 &=&  \frac{12109\omega^2+31536}{1393056}y^2 +\frac{107\omega^2+2106}{696528}Y^2  \\ \nonumber
&& +\frac{113\omega^2+918}{29232}z^2 +\frac{4\omega^2+81}{14616}Z^2,
\\ \label{deltaxx}
\Delta_2 &=& \frac{9}{38696}(13 \omega^2+159)y Y +\frac{3}{1624}(3\omega^2+10)z Z,
\end{eqnarray}
whose right members must be evaluated using the variables with subindex 2 for the direct corrections in Eq.~(\ref{directCart}), and using the variables with subindex 1 for the inverse corrections $x_2=x_1+\rho^2\sigma(\lambda\Delta_1+\Delta_2)$, $X_2=X_1-\rho^2\sigma(\lambda\Delta_1-\Delta_2)$. 
\par

\subsection{Detuning and Lissajous variables}

The unperturbed frequency of the oscillations in the $z$ direction can be written as $\nu=\omega\sqrt{1-\delta}$ where, in view of Eqs.~(\ref{nu}) and (\ref{lw}),
\[
\delta=1-(\nu/\omega)^2=\frac{23-8\sqrt{7}}{27}\approx0.068,
\]
is a ``detuning'' parameter \citep{Henrard1970} that amounts to one tenth of the Hill radius $\rho=3^{-1/3}$, and will be taken as a first order perturbation. Then, the Hamiltonian of the center manifold in Eq.~(\ref{centerHam}) is rearranged in the form
\begin{equation} \label{Hampeo}
\mathcal{C}=\mbox{$\frac{1}{2}$}(Y_2^2+Z_2^2)+\mbox{$\frac{1}{2}$}\omega^2(y_2^2+z_2^2)+\tilde{C}_1+\frac{1}{2!}C_2, 
\end{equation}
with $\tilde{C}_1\equiv{C}_1-\frac{1}{2}\omega^2\delta z_2^2$.
\par

Equation (\ref{Hampeo}) can be viewed as the Hamiltonian of a perturbed elliptic oscillator whose principal part comprises two harmonic oscillators in the 1-1 resonance. Moreover, because the perturbation belongs to the real algebra in the Cartesian variables $(y_2,z_2,Y_2,Z_2)$, the Hamiltonian (\ref{Hampeo}) is advantageously attacked in Lissajous variables \citep{Deprit1991}.
\par

The Lissajous transformation
\begin{equation} \label{T4}
\mathcal{T}_4:(y_2,z_2,Y_2,Z_2)\longrightarrow(\ell,g,L,G;\omega),
\end{equation}
is defined as
\begin{eqnarray} \label{yLi}
y_2 &=& s \cos (g+\ell )-d \cos (g-\ell ), \\ \label{zLi} 
z_2 &=& s \sin (g+\ell )-d \sin (g-\ell ), \\ \label{yyLi} 
Y_2 &=&-\omega\left[s \sin (g+\ell )+d \sin (g-\ell )\right], \\ \label{zzLi} 
Z_2 &=&\phantom{-}\omega\left[s \cos (g+\ell )+d \cos (g-\ell )\right],
\end{eqnarray}
where $s\equiv{s}(L,G;\omega)$ and $d\equiv{d}(L,G;\omega)$ are the state functions
\begin{equation} \label{sdLi}
s=\sqrt{\frac{L+G}{2\omega}}, \qquad d=\sqrt{\frac{L-G}{2\omega}}.
\end{equation}
\par

The variables in Eqs.~(\ref{yLi})--(\ref{sdLi}) have full geometrical meaning: They define an ellipse in the $y_2$-$z_2$ plane centered at the origin, whose size and shape are defined by the semi-major axis $a$ and semi-minor axis $|b|$ that are derived from the relations $L=\frac{1}{2}\omega(a^2+b^2)$, $G=\omega{a}b$, with the direction of the semi-minor axis with respect to the axis of ordinates defined by the angle $g$; in this ellipse, the elliptic anomaly $\ell$ is measured from the semi-major axis $b$. An analogous ellipse is defined also in the $Y_2$-$Z_2$ plane, now with semi-major axis $\omega{a}$ and semi-minor one $\omega|b|$ (see p.~209 and ff.~of \citet{Deprit1991} for full details).
\par

The Lissajous transformation is applied to Eq.~(\ref{Hampeo}), to give 
\[
\mathcal{T}_4\circ\mathcal{C}(y_2,z_2,Y_2,Z_2)=\mathcal{A}(\ell,g,L,G)\equiv\sum_{n=0}^{2}\frac{\epsilon^n}{n!}\mathcal{A}_n,
\]
where
\begin{eqnarray}
\mathcal{A}_0 &=& \omega L \\ \label{L1}
\mathcal{A}_1 &=& \frac{1}{4}\delta\omega^2\sum_{i=0}^{1}\sum_{j=-1}^{1}Q_{1,2i,2j}\cos(2ig+2j\ell) \\ \nonumber
&& +\frac{3}{448} \rho ^2 \tau  \omega\sum_{i=0}^{1}\sum_{j=-1}^{2}Q_{1,2i+1,2j-1}\sin[(2i+1)g+(2j-1)\ell] \\ \label{L2}
\mathcal{A}_2 &=& \frac{9 \rho }{62842304}\sum_{i=0}^{2}\sum_{j=-2}^{2}Q_{2,i,j}\cos(2ig+2j\ell) 
\end{eqnarray}
and the coefficients $Q_{n,j,k}$, which only depend on the momenta $L$ and $G$ through the state functions $s$ and $d$, are given in Table \ref{t:Lo23}.

\begin{table}[htb]
\begin{tabular}{@{}ll@{}}
\hline
$Q_{1,0,-2}=d s$ & $Q_{1,2,-2}=d^2$ \hspace{2.8cm}\vphantom{$\frac{M^N}{1}$}\\
$Q_{1,0,0}=-d^2-s^2$ & $Q_{1,2,0}=-2 d s$ \\
$Q_{1,0,2}=d s$ & $Q_{1,2,2}=s^2$ 
\\[1ex]
$Q_{1,1,-3}=(7-10\omega^2)d^2 s$ & $Q_{1,3,-3}=(11-2\omega^2)d^3$ \\
$Q_{1,1,-1}=(86-20\omega^2)ds^2+(6\omega^2+3)d^3$ & $Q_{1,3,-1}=(10\omega^2-43)d^2s$ \\
$Q_{1,1,1}=(86-20\omega^2)d^2s+(6\omega^2+3)s^3$ & $Q_{1,3,1}=(10\omega^2-43)ds^2$ \\
$Q_{1,1,3}=(7-10\omega^2)ds^2$ & $Q_{1,3,3}=(11-2\omega^2)s^3$ \\
\end{tabular}
\begin{tabular}{@{}ll@{}}
\hline
$Q_{2,0,\pm2}=\frac{9}{4}(454826\omega^2-29767905)d^2s^2$ \vphantom{$\frac{M^N}{1}$} \\[0.33ex]
$Q_{2,0,\pm1}=\frac{3}{2}(4601090\omega^2-7248069)(d^2+s^2)ds$ \\ [0.33ex]
$Q_{2,0,0}=\frac{3}{4}(54449757-4733570\omega^2)(d^4+s^4)-27(7866699-1473422\omega^2)d^2s^2$ \\ [0.9ex]
$Q_{2,1,\pm2}=(4343546\omega^2+13096395)d^{2\mp1}s^{2\pm1}$ \\ [0.33ex]
$Q_{2,1,\pm1}=3(56290797-12398698\omega^2)d^2s^2-(45973827-8226998\omega^2)d^{2\mp2}s^{2\pm2}$ \\ [0.33ex]
$Q_{2,1,0}=-9(1640482\omega^2-5859465)(d^2+s^2)ds$ \\ [0.9ex]
$Q_{2,2,\pm2}=\frac{1}{4}(35821341-9394466\omega^2)d^{2\mp2}s^{2\pm2}$ \\ [0.33ex]
$Q_{2,2,\pm1}=(3324843-622142\omega^2)d^{2\mp1}s^{2\pm1}$ \\ [0.33ex]
$Q_{2,2,0}=\frac{27}{2}(1658222\omega^2-6951883)d^{2}s^{2}$ \\ [0.33ex]
\hline
\end{tabular}
\caption{Coefficients $Q_{n,j,k}$ in Eqs.~(\protect\ref{L1}) and (\protect\ref{L2}) \label{t:Lo23}. }
\end{table}

\subsection{Elimination of the eccentric anomaly}

The Hamiltonian in Lissajous variables can be reduced to a one degree of freedom Hamiltonian by means of a new canonical transformation
\begin{equation} \label{T5}
\mathcal{T}_5:(\ell,g,L,G)\longrightarrow(\ell',g',L',G';\epsilon)
\end{equation}
such that, after truncation to $\mathcal{O}(\epsilon^2)$,
\begin{equation} \label{HamB}
\mathcal{T}_5\circ\mathcal{A}(\ell,g,L,G)=\mathcal{B}(-,g',L',G')\equiv\sum_{n=0}^{2}\frac{\epsilon^n}{n!}B_n,
\end{equation}
where,
\begin{eqnarray} \label{B0}
B_0 &=& \omega{L}', \\ \label{B1}
B_1 &=& -\frac{1}{4}\delta\omega(L'+2\omega{d}'s'\cos2g'), \\ \label{B2}
B_2 &=& 2!\left(\frac{\delta}{4}B_1
-k_1L'^2 +k_2L'\omega{s}'d'\cos2g'-k_3\omega^2s'^2d'^2\cos4g' +\frac{k_4}{4}G'^2\right),\qquad
\end{eqnarray}
where primes in functions mean the same functions written in the prime variables, and
\begin{eqnarray*}
k_1 &=& \frac{1}{16} (6829135-609646\omega^2)k_0, \\
k_2 &=& (126184-9583\omega^2)k_0, \\
k_3 &=&-\frac{3}{4} (439957-103954\omega^2)k_0, \\
k_4 &=& \frac{3}{4} (7293079-1280862\omega^2)k_0, \\
k_0 &=& \frac{1}{6733104}(\omega^2+2)\rho,
\end{eqnarray*}
are strictly positive irrational numbers ($k_1\approx0.17$, $k_2\approx0.055$, $k_3\approx0.003$, $k_4\approx0.87$, $k_0\approx6.5\times10^{-7}$) that have been introduced for abbreviating expressions.
\par

Up to the first order, the short-period corrections of the transformation (\ref{T5}) are of the form
\begin{eqnarray} \label{shortperiod}
\Delta\xi &=& \frac{\delta}{4}\sum_{i=-1}^{1}\frac{\xi_{2,i}}{4s'd'}\cos(2ig'+2\ell'-\beta) \\ \nonumber
&& +\frac{\rho^2\tau\omega}{4032}\sum_{i=0}^{1}\sum_{j=-1}^{2}\frac{\xi_{2i+1,2j-1}}{4s'd'}\sin[(2i+1)g'+(2j-1)\ell'+\beta],
\end{eqnarray}
where $\xi\in(\ell',g',L',G')$, and $\beta=0$ for the momenta while $\beta=\pi/2$ in the case of the coordinates. The necessary coefficients are given in Table \ref{t:dlgLG1}, where primes have been dropped for brevity. Note that $G_{i,j}=(m/n)L_{i,j}$ where $m$ is the coefficient of $g'$ and $n$ is the coefficient of $\ell'$ in the argument of the trigonometric function factored by $G_{i,j}$, and hence the corresponding coefficients are not provided.
\par

\begin{table}[htb]
\begin{tabular}{@{}lllll@{}}
$i,j$ & $L_{i,j}/(4sd\omega)$ & $\ell_{i,j}$ & $g_{i,j}$ & $c_{i,j}$
\\ [0.33ex]
\hline
${1,-3}$ & $-3c_{1,3}d^2s$ & $c_{1,3}(d^2+2s^2)d$ & $c_{1,3}(d^2-2s^2)d$ \vphantom{$\frac{M^N}{1}$} \\ 
${1,-1}$ & $3c_0d^3-2c_{3,1}ds^2$ & $2c_{3,1}s^3-c_2sd^2$ & $8c_{3,1}d^2s-\ell _{1,-1}$ \\
${1,1}$ & $3c_0s^3-2c_{3,1}d^2s$ & $c_2ds^2-2c_{3,1}d^3$ & $8c_{3,1}ds^2+\ell _{1,1}$ \\ 
${1,3}$ & $-3c_{1,3}ds^2$ & $-c_{1,3}(2d^2+s^2)s$ & $c_{1,3}(s^2-2d^2)s$ & $\frac{7}{3}\omega ^2-\frac{256}{3}$ \\ 
$2,-1$ & $-d^2$ & $d s$ & $-d s$ \\
${2,0}$ & $-2 d s$ & $d^2+s^2$ & $d^2-s^2$ \\ 
${2,1}$ & $-s^2$ & $d s$ & $d s$ \\ 
${3,-3}$ & $-c_{3,3}d^3$ & $c_{3,3}d^2s$ & $-c_{3,3}d^2s$ \\ 
${3,-1}$ & $c_{3,1}d^2s$ & $-c_{3,1}(d^2+2s^2)d$ & $c_{3,1}(2s^2-d^2)d$ \\ 
${3,1}$ & $c_{3,1}ds^2$ & $c_{3,1}(2d^2+s^2)s$ & $c_{3,1}(2d^2-s^2)s$ & $43 \omega ^2-184$\\ 
${3,3}$ & $-c_{3,3}s^3$ & $-c_{3,3}ds^2$ & $-c_{3,3}ds^2$ & $11 \omega ^2-32$\\
\hline
\end{tabular}
\caption{Coefficients $\xi_{i,j}$ in Eqs.~(\protect\ref{shortperiod}); $c_0=c_{3,1}-4 c_{3,3}$ and $c_2=5c_{3,1}-36c_{3,3}$. \label{t:dlgLG1} }
\end{table}

\section{The reduced phase space} \label{se:reduflow}

As a result of the averaging the elliptic anomaly beocmes cyclic in Eq.~(\ref{HamB}), and, therefore, its conjugate momentum $L'$ is an integral of the motion. Thus, the problem has been reduced to a one degree of freedom Hamiltonian in the coordinate $g'$ and its conjugate momentum $G'$.
\par

The $(g,G)$ chart is a cylindrical map that misses the circular orbits representation. Indeed, when the ellipse's semi-major and semi-minor axes are equal the angle $g$ is undetermined and the Lissajous transformation is singular. This fact does not invalidates the normalization that has been carried out to eliminate $\ell$. Since the perturbation method used is invariant with respect to canonical transformations \citep{Hori1966}, the resulting theory remains valid after reformulated in nonsingular variables, which, besides, do not need to be canonical, cf.~\citep{DepritRom1970}.
\par

Therefore, circular orbits are not excluded from the theory, and the reduced phase space is conveniently analyzed in terms of the Hopf coordinates. Thus, the new transformation
\begin{equation}
\mathcal{T}_6:(g',G';L';\omega)\longrightarrow(I_1,I_2,I_3),
\end{equation}
defined by
\begin{equation} \label{I123}
I_1 = \omega s' d' \cos2g', \qquad
I_2 = \omega s' d' \sin2g', \qquad
I_3 = \mbox{$\frac{1}{2}$} G',
\end{equation}
which projects the reduced phase space onto the sphere centered at the origin
\begin{equation} \label{sphere}
I_1^2+I_2^2+I_3^2=\mbox{$\frac{1}{4}$}L'^2,
\end{equation}
is applied to the reduced Hamiltonian $\mathcal{B}$, yielding
\begin{equation} \label{HamHopf}
\mathcal{I} = \omega(1-\mbox{$\frac{1}{2}$}\delta^*)L'-k_1L'^2 +(k_2L'-\omega\delta^*)I_1+k_3(I_2^2-I_1^2)+k_4I_3^2,
\end{equation}
where the new abbreviation
\[
\delta^*=\mbox{$\frac{1}{2}$}\delta\left(1+\mbox{$\frac{1}{4}$}\delta\right)\approx0.03454,
\]
has been introduced.
\par

In the general case, the Hamiltonian flow associated to Eq.~(\ref{HamHopf}) is given by the differential system
\begin{eqnarray} \label{I1p}
\dot{I}_1 &=& 2(k_3-k_4)I_2I_3, \\ \label{I2p}
\dot{I}_2 &=& \left[\delta^*\omega -k_2L'+2(k_3+k_4){I}_1\right]{I}_3, \\ \label{I3p}
\dot{I}_3 &=&-\left(\delta^*\omega -k_2L'+4k_3{I}_1\right){I}_2.
\end{eqnarray}
The particular case $L'=\tilde{L}=\delta^*\omega/k_2\approx1.29839$ yields a differential system that is analogous to the Euler equations for the free rigid body motion, and, therefore, can be integrated in terms of Jacobi elliptic functions. The general solution of Eqs.~(\ref{I1p})--(\ref{I3p}) is not pursued. If found, it would not provide much insight into the nature of the solution due to the unavoidable use of special functions. However, a lot of qualitative and quantitative information can be obtained from the study of particular solutions, like the equilibria, as well as from the graphic representation of the flow.
\par

\subsection{Equilibria}

When $I_2=I_3=0$ Eqs.~(\ref{I1p})--(\ref{I3p}) vanish. In consequence, the points of the sphere
\begin{equation} \label{Epm1}
E_{\pm1}=\left(\pm\mbox{$\frac{1}{2}$}L',0,0\right),
\end{equation}
are always equilibria of the reduced system. On the other hand, when $I_2=0$ but $I_3\ne0$, two new equilibria
\begin{equation} \label{Epm2}
E_{\pm2}=\left(I_{1,\mathrm{H}},0,\pm\mbox{$\frac{1}{2}$}(L'^2-4I_{1,\mathrm{H}}^2)^{1/2}\right),
\end{equation}
may exist, where $I_{1,\mathrm{H}}$ is the solution of $\dot{I}_2=0$ with $I_3\ne0$, namely
\begin{equation} \label{I1H}
I_{1,\mathrm{H}} =\frac{k_2L'-\delta^*\omega}{2(k_3+k_4)},
\end{equation}
which vanishes for the particular value $\tilde{L}$ mentioned before, in which case $I_3=\frac{1}{2}\tilde{L}$. In addition, $I_{1,\mathrm{H}}$ needs to fulfill the geometric condition $|I_1|\le\frac{1}{2}L'$, from which
\begin{equation} \label{L0halo}
L'\ge{L}_0=\frac{\delta^*\omega}{k_2+k_3+k_4}\approx0.0768606.
\end{equation}
By replacing $L_0$ into Eq.~(\ref{I1H}) it is shown that the new equilibria bifurcate from $E_{-1}$ when $I_{1,\mathrm{H}}=-\frac{1}{2}L_0\approx-0.0384303$.
\par

Analogously, when $I_3=0$, $I_2\ne0$, two new equilibria
\begin{equation} \label{Epm3}
E_{\pm3}=\left(I_{1,\mathrm{B}},\pm\mbox{$\frac{1}{2}$}(L'^2-4I_{1,\mathrm{B}}^2)^{1/2},0\right),
\end{equation}
may exist, where $I_{1,\mathrm{B}}$ is the solution of $\dot{I}_3=0$ with $I_2\ne0$, viz.
\begin{equation} \label{I1B}
I_{1,\mathrm{B}} =\frac{k_2L'-\delta^*\omega}{4k_3},
\end{equation}
which, again, vanishes for $L'=\tilde{L}$, a case in which $I_2=\frac{1}{2}\tilde{L}$. Note that $(k_3+k_4)I_{1,\mathrm{H}}=2k_3I_{1,\mathrm{B}}$. Again, the geometry of the sphere introduces the constraint
\begin{equation} \label{L12bridge}
\frac{\delta^*\omega}{k_2-2k_3}=L_1\le L' \le L_2=\frac{\delta^*\omega}{k_2+2k_3}.
\end{equation}
By replacing $L'=L_1\approx1.17113$ and $L'=L_2\approx1.45668$ into Eq.~(\ref{I1B}), the sign taken by $I_{1,\mathrm{B}}$ shows that the bifurcation at $L_1$ happens from $E_{-1}$, whereas at $L_2$ the bifurcation occurs from $E_{+1}$. Thus, the symmetric equilibria $E_{\pm3}$ migrate from $E_{-1}$ to $E_{+1}$ from increasing values of $L'$, or vice-versa when $L'$ decreases from $L_2$ to $L_1$.
\par

The stability of the equilibria can be computed from the usual linearization of the flow in Eqs.~(\ref{I1p})--(\ref{I3p}). It shows that $E_1$ is of the elliptic type for $L'<L_2$, and then of the hyperbolic type. Besides, $E_{-1}$ is of the hyperbolic type between $L_0$ and $L_1$, and elliptic otherwise. $E_{\pm2}$ are stable when they exist ($L'>L_0$), whereas $E_{\pm3}$ are unstable from their bifurcation from $E_{-1}$ at $L'=L_1$ until they collapse with $E_{1}$ at $L'=L_2$. This behavior will be confirmed in the next section visualizing the flow. 
\par

Besides, it is noted that for small enough values of $L'$, and, in consequence of $I_1$, $I_2$, and $I_3$, effects of $\mathcal{O}(L'^2)$ may be neglected. Then the first order truncation of the Hamiltonian (\ref{HamB}) can be taken as representative of the motion, namely $\mathcal{B}\approx{B}_0+B_1$, which, in the Hopf variables is
$\mathcal{I}\approx\omega(1-\frac{1}{4}\delta){L}'-\frac{1}{2}\delta\omega I_1$.
Then, $\dot{I}_1\approx0$, $\dot{I}_2\approx\frac{1}{2}\delta\omega{I}_3$, $\dot{I}_3\approx-\frac{1}{2}\delta\omega{I}_2$, which show that $I_1$ remains approximately constant along each trajectory on the sphere, whereas $I_2$ and $I_3$ evolve with the slow frequency $\Omega=\frac{1}{2}\delta\omega$ in circumferences parallel to the plane $I_1=0$ of radius $\frac{1}{2}(L'^2-4I_1^2)^{1/2}$. 
\par

A caveat is in order at this point. While the analytical solution has been constrained to the lower orders in the expansion of the Hamiltonian at the libration point, and, in consequence, to compatible orders of the perturbation theories carried out, the bifurcations predicted by the theory occur at high values of $L'$. Hence, this information must be taken just as qualitative and keep in mind that the computed bifurcation values may appreciably change when using a higher order theory. 

\subsection{Visualizing the flow}

On the other hand, the changes in the reduced flow discussed previously can be easily visualized in the sphere without need of integrating Eqs.~(\ref{I1p})--(\ref{I3p}). Indeed, for a given value of the dynamical parameter $L'$, a trajectory in the manifold $\mathcal{I}=h$ is defined by the intersection of the two-dimensional surface defined by the Hamiltonian (\ref{HamHopf}) with the surface defined by the sphere (\ref{sphere}). Hence, by elimination of $I_3$ between these two equations,
\begin{equation} \label{I22}
I_2^2=-\frac{(k_3+k_4)I_1^2+(\delta^*\omega-k_2L')I_1+(k_1-\mbox{$\frac{1}{4}$}k_4)L'^2-(1-\mbox{$\frac{1}{2}$}\delta^*)\omega{L}'+h}{k_4-k_3},
\end{equation}
which provides $I_2$ as a function of $I_1$ and the pair of dynamical parameters $L'$ and, $h$. An analogous elimination of $I_2$ yields
\begin{equation} \label{I32}
I_3^2=\frac{2k_3I_1^2+(\delta^*\omega-k_2L')I_1+(k_1-\mbox{$\frac{1}{4}$}k_3)L'^2-(1-\mbox{$\frac{1}{2}$}\delta^*)\omega{L}'+h}{k_4-k_3},
\end{equation}
which gives $I_3\equiv{I}_3(I_1;L',h)$. Therefore, each trajectory $\mathcal{I}(I_1,I_2,I_3;L')=h$ can be depicted on the sphere from the simple evaluation of the square roots of Eqs.~(\ref{I22}) and (\ref{I32}) in  such subset of the interval $I_1\in[-\frac{1}{2}L',\frac{1}{2}L']$ in which the square roots are real.
\par

\begin{figure}
\centering \includegraphics[scale=0.7, angle=-90]{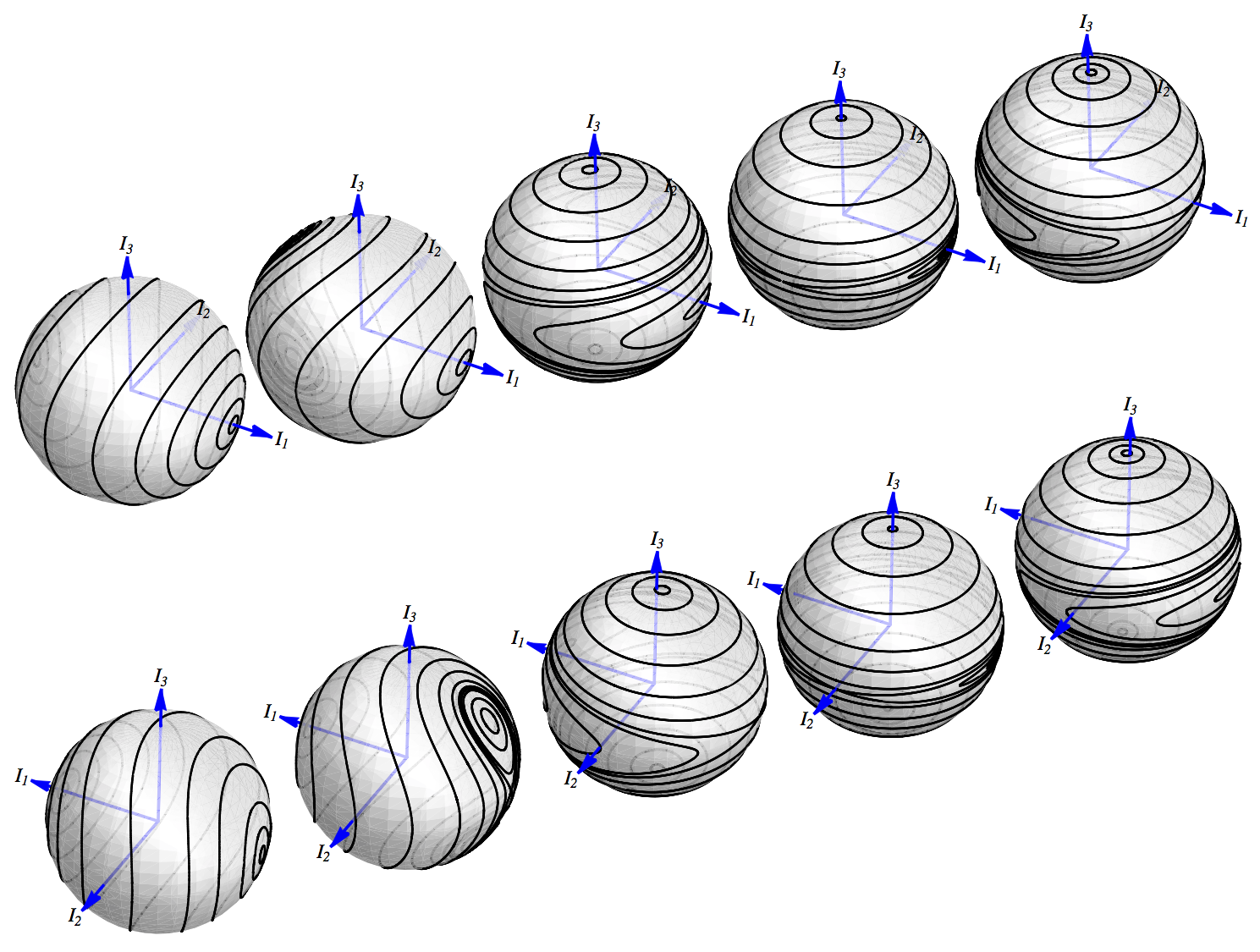}
\caption{Opposite views (left and right columns) of the bifurcation sequence of the Hamiltonian flow in Eq.~(\protect\ref{HamHopf}) for increasing values of $L'$. From top to bottom, $L'=0.05,0.1,0.7,1.29839$ and $2.5$.}
\label{f:spheresSequence}
\end{figure}

The sequence of bifurcations of the flow of the reduced problem presented in Fig.~\ref{f:spheresSequence} has been depicted using this technique. Two different views of each sphere are shown in the figure, the second one obtained by rotating the first one 180 degrees about the axis $I_3$. In spite of the various pairs of spheres correspond to different values of the dynamical parameter $L'$, and, in consequence, should have different radius, they are represented with a normalized radius 1 to better appreciate the flow. As shown in Fig.~\ref{f:spheresSequence}, for small values of $L'$ the two equilibria $E_{\pm1}$ are stable and the flow circulates about them (first row of Fig.~\ref{f:spheresSequence}). For increasing values of $L'$, the flow distorts about $E_{-1}$, until, eventually, $L_{-1}$ changes to instability in a bifurcation event, and two new equilibria appear in the plane $I_2=0$ (second row of Fig.~\ref{f:spheresSequence}). The bifurcated equilibria move on the $I_1$-$I_3$ meridian towards the $\pm{I}_{3}$ axis, while the flow narrows about the plane $I_3=0$ plane (third row of Fig.~\ref{f:spheresSequence}). Eventually, the equilibrium $E_{-1}$ comes back to stability in a new bifurcation, and two new unstable equilibria appear in the plane $I_3=0$ (fourth row of Fig.~\ref{f:spheresSequence}), which migrate along the equator of the sphere, the $I_1$-$I_2$ circumference, for increasing values of $L'$ until merging in a new bifurcation phenomenon with the equilibrium $E_{+1}$, which undergoes a concomitant change to instability (fifth row of Fig.~\ref{f:spheresSequence}). Further increases of $L'$ do not introduce qualitative changes in the flow. 

\subsection{Orbits of the center manifold}

The orbits of the reduced phase space are now identified with orbits of the original Hill problem. From Eqs.~(\ref{Epm1}), (\ref{Epm2}), and (\ref{Epm3}), and the inverse transformation of Eq.~(\ref{I123}) one finds that $G'=0$ for $E_{\pm1}$ and $E_{\pm3}$. Hence, except for the short-period effects due to the corrections in Eq.~(\ref{shortperiod}), the corresponding motion in the center manifold is rectilinear. On the contrary, $G'\ne0$ for $E_{\pm2}$, which, therefore, correspond to elliptic motion, on average, in the center manifold. Besides, $g'=0$ for $E_{1}$, which in view of Eqs.~(\ref{yLi})--(\ref{zzLi}), yields, on average, harmonic oscillations in the $z_2$ direction, whereas $g'=\frac{\pi}{2}$ for $E_{-1}$ thus constraining the oscillations to the $y_2$ direction. On the other hand, the solutions $E_{\pm2}$ yield elliptic oscillations in the $y_2$-$z_2$ plane, the area of the respective ellipses depending on $|G'|=(L'^2-4I_{1,\mathrm{H}}^2)^{1/2}$. Finally, the equilibria $E_{\pm3}$ result in rectilinear oscillations in the $y_2$-$z_2$ plane with inclination $g'$ given by the components $I_1$ and $I_2$ of the equilibria.
\par

When the Cartesian coordinates are recovered through the linear transformation in Eq.~(\ref{x1tox}), one recognizes that the the equilibria of the reduced problem correspond to the well known periodic orbits of the center manifold of the Hill problem \citep{GomezMarcoteMondelo2005}. Namely,
\begin{itemize}
\item vertical ($E_{1}$) and planar ($E_{-1}$) Lyapunov orbits
\item Halo orbits, with the upper part towards the libration point ($E_{2}$) and the symmetric one ($E_{-2}$) with the upper part towards the primary
\item orbits of the two lane bridge ($E_{\pm3}$) connecting vertical and planar Lyapunov orbits
\end{itemize}
\par

The periods of these periodic orbits is estimated from the rate of variation of the eccentric anomaly in the Lissajous normalized coordinates $T=2\pi/\dot{\ell}'$, where $\dot\ell'=\partial\mathcal{B}/\partial{L}'$ is computed from Eq.~(\ref{HamB}). Up to the first order in the small parameter
\begin{equation}
\dot\ell=\omega-\omega\frac{1}{4}\delta\left(1+\frac{1}{\sqrt{1-G'^2/L'^2}}\cos2g'\right),
\end{equation}
where in each case the values of $g'$, $G'$ and $L'$ must be replaced from those of the corresponding equilibria.
\par

\section{Validation tests} \label{se:validation}

Some examples are provided to illustrate the application and performance of the analytical theory. 

\begin{figure}[htb]
\centering 
\includegraphics[scale=0.71, angle=0]{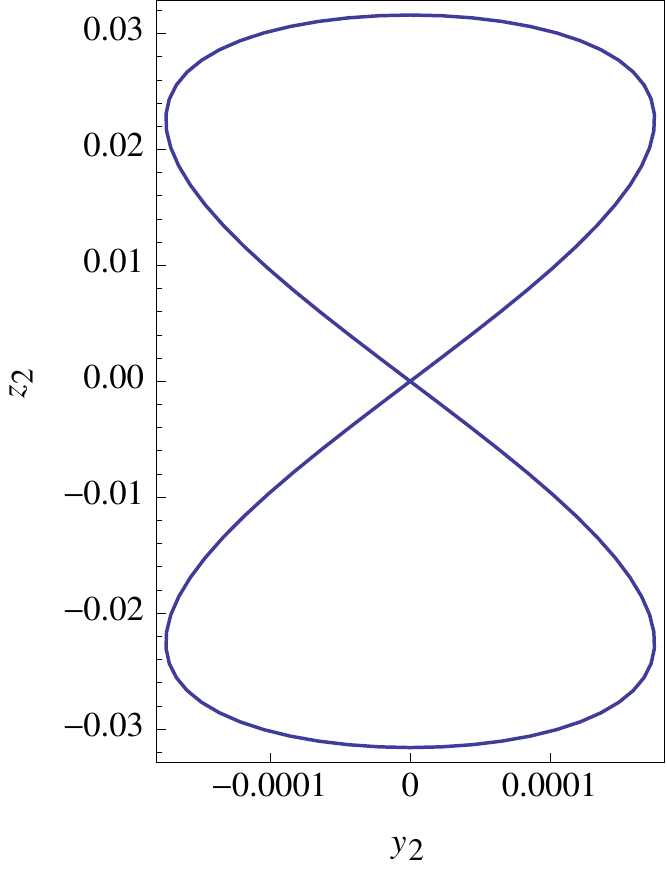}\qquad
\includegraphics[scale=0.71, angle=0]{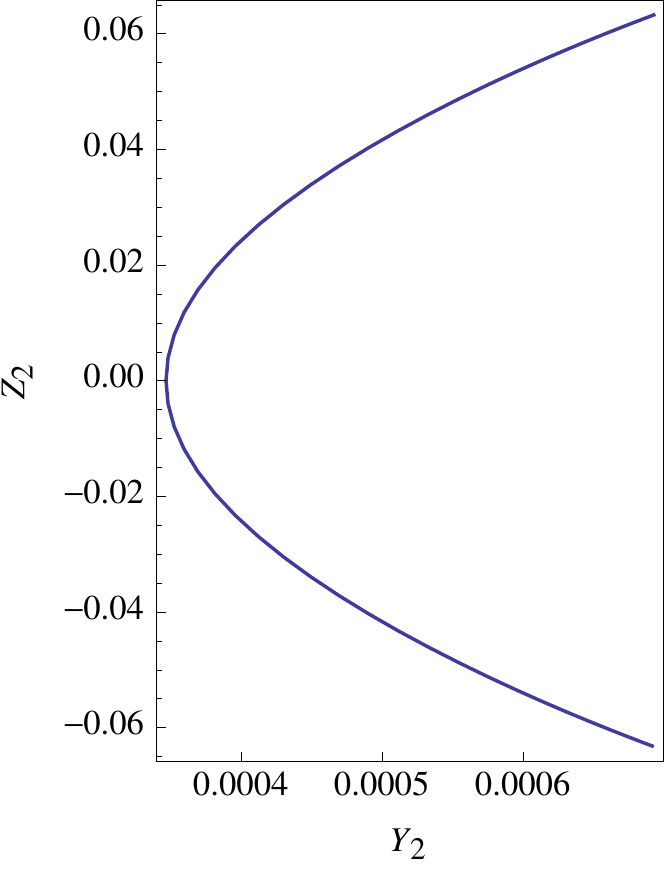}
\caption{Lyapunov vertical solution for $L=0.001$ in the center manifold. Left: coordinates space. Right: momenta space. Note the different scales in abscissas and ordinates.}
\label{f:LyVsmall}
\end{figure}

Thus, starting from $L'=0.001$, the Hopf coordinates of the $E_{1}$ equilibrium, corresponding to the Lyapunov vertical solution, are computed from Eq.~(\ref{Epm1}): $I_1=0.0005$, $I_2=I_3=0$, what result into $g'=0$ and $G'=0$, to which correspond a period $T=3.13965$ of the rectilinear oscillations. Then, for each value $\ell'\in[0,2\pi)$, the original Lissajous variables are computed by recovering the short-period corrections in Eq.~(\ref{shortperiod}). Then, the Lissajous transformation in Eqs.~(\ref{yLi})--(\ref{zzLi}) provides the coordinates $y_2$, $z_2$, $Y_2$, $Z_2$, in the center manifold. Figure~\ref{f:LyVsmall} illustrates how this orbit looks like in the center manifold.
\par

\begin{figure}[htb]
\centering 
\includegraphics[scale=0.8, angle=0]{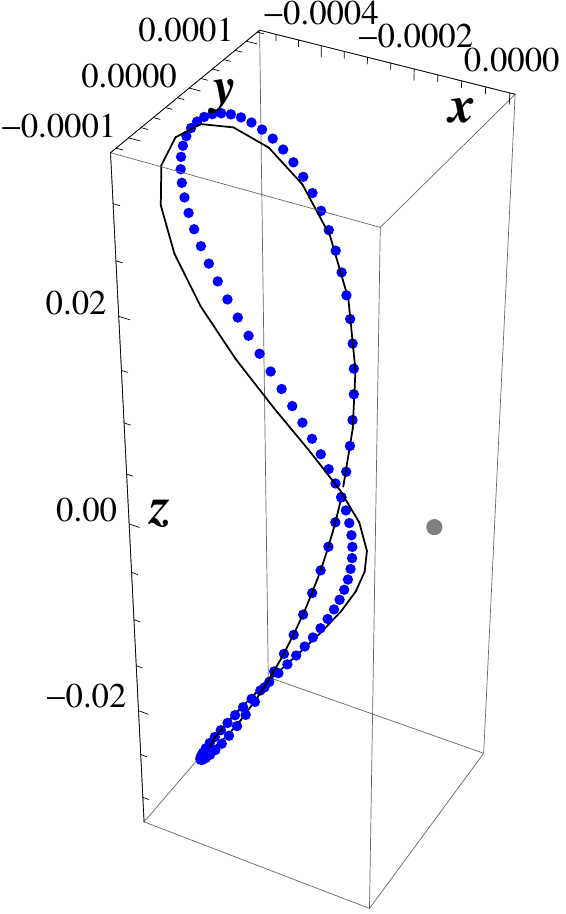} \qquad
\includegraphics[scale=0.77, angle=0]{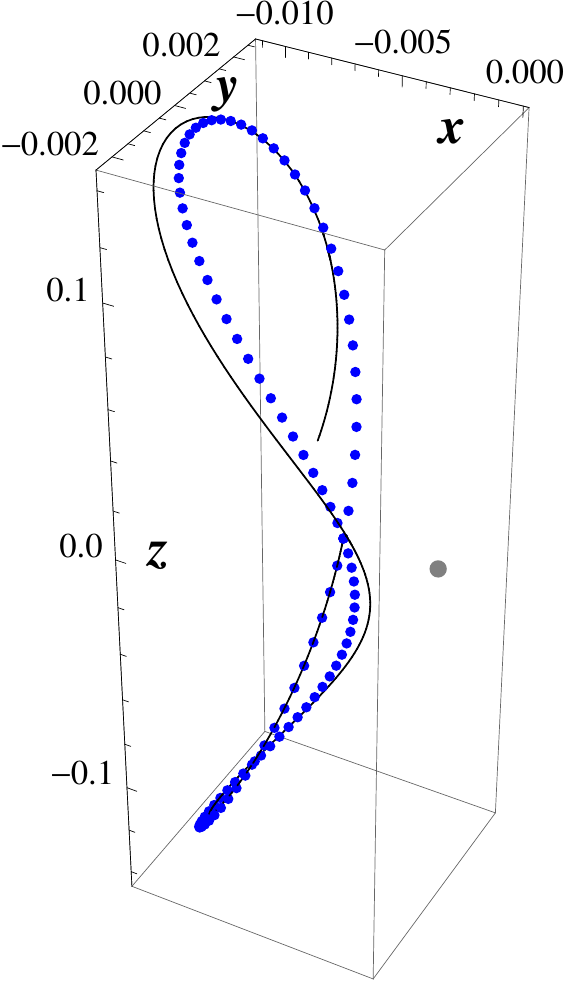}
\caption{Lyapunov vertical orbits in Cartesian coordinates. Left: $L=0.001$, right: $L=0.02$. The orbits are distorted by the scaling of the $z$ axis. The origin is the libration point, which is highlighted with a gray dot.}
\label{f:vLy001orbit}
\end{figure}

To recover the orbit in the original space, the corrections in Eqs.~(\ref{deltax}) and (\ref{deltaxx}) are computed first to get the subindex 1 variables, and then Eq.~(\ref{x1tox}) provides the Cartesian coordinates of the orbit in the original space, which is illustrated in the left plot of Fig.~\ref{f:vLy001orbit}. The analytical solution obtained in this way is obviously periodic by construction (blue points in the left plot of Fig.~\ref{f:vLy001orbit}). However, when the initial conditions obtained from the analytical solution for, say, $\ell'=0$ are propagated in the original equations of motion, Eqs.~(\ref{dpx})--(\ref{dPz}), the orbit is not exactly periodic due to the neglected higher order corrections in the transformations computed by perturbation theory (black curve in the left plot of Fig.~\ref{f:vLy001orbit}), namely, the reduction to the center manifold and the averaging of the elliptic anomaly, as well as in the computation of the period. Indeed, the differences between the initial state and the computed state after the period given by the analytical approximation is of the order of $10^{-5}$ in this example. 
\par

As expected, because the analytical solution is constrained to small values of the distance when compared to the Hill radius, this lack of periodicity is more evident for higher values of $L'$, because of the corresponding larger size of the orbit. The differences between the analytical solution and the true orbit for the same initial conditions soon become evident, as illustrated in the right plot of Fig.~\ref{f:vLy001orbit} for $L'=0.02$, where dots correspond to the analytical orbit and the black curve represents to the numerical solution. Still, for moderate values of $L'$ the analytical solution provides good enough seeds to feed a differential corrections procedure that easily gets the periodic solution. Thus, even though the periodicity error when using initial conditions of the analytical theory is of the order of $10^{-2}$ for $L'=0.02$, the differential corrections algorithm in \citep{LaraPelaez2002} only needs to compute four consecutive corrections to converge to a periodic orbit with a periodicity error of $\mathcal{O}(10^{-13})$.
\par

\begin{figure}[htb]
\centering 
\includegraphics[scale=0.71, angle=0]{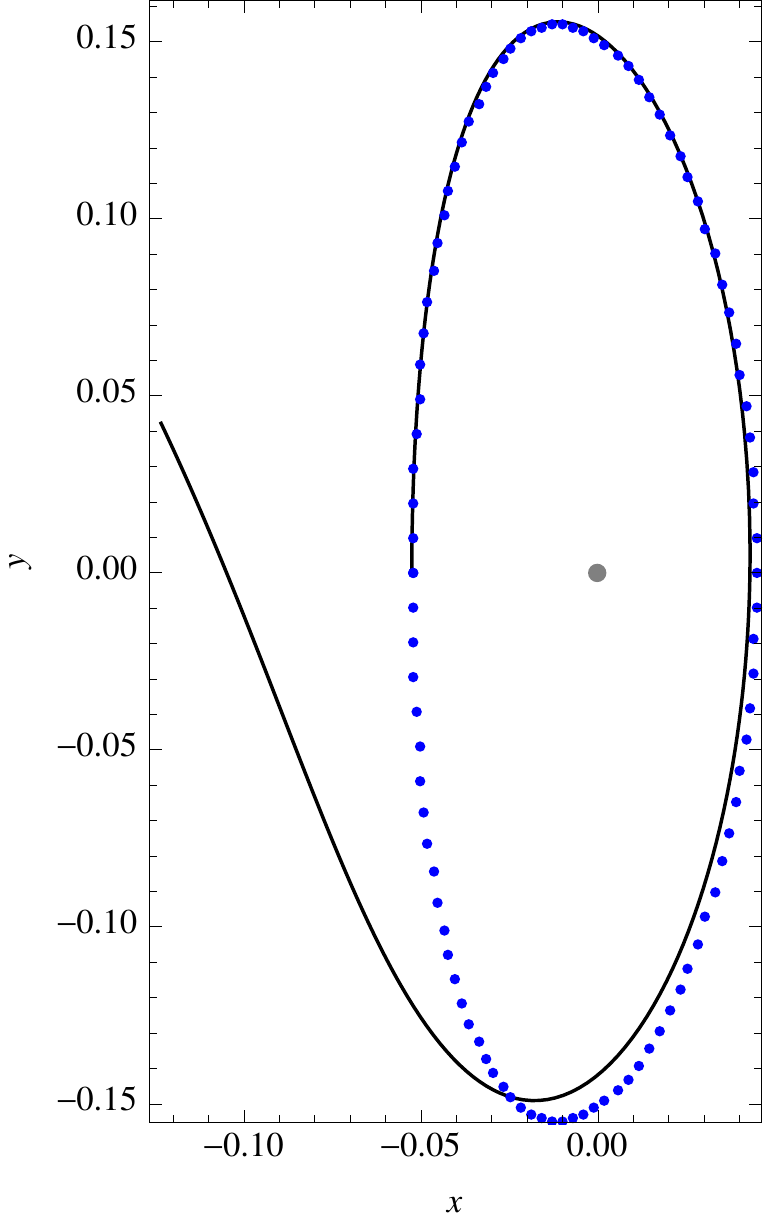}\qquad
\includegraphics[scale=0.7, angle=0]{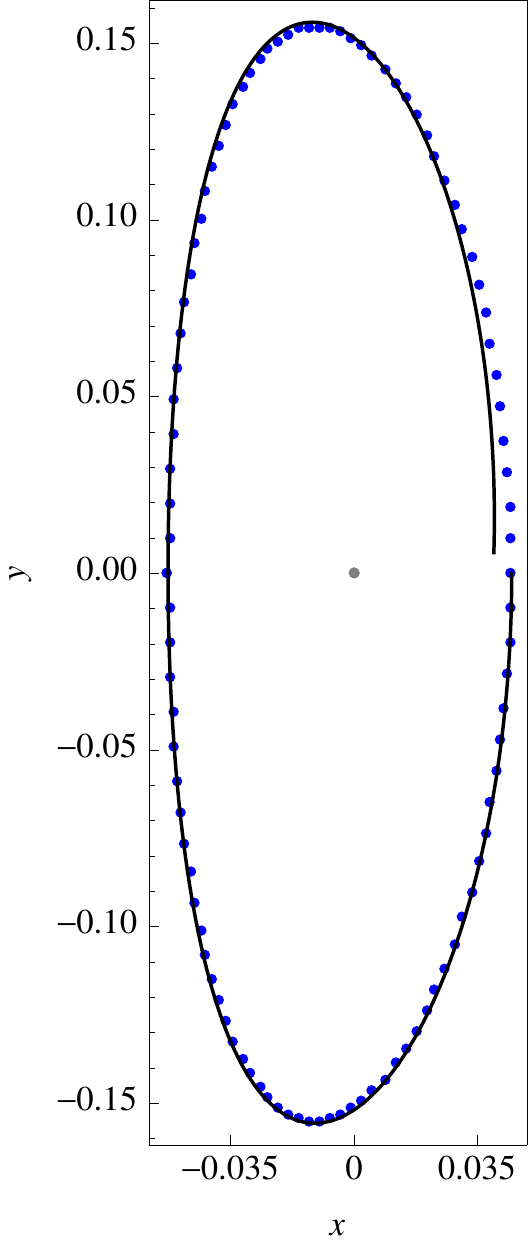}
\caption{Planar Lyapunov orbits in Cartesian coordinates for $L=0.02$ and different initial conditions. Blue points: analytical solution; black curves: numerical propagation of the Hill problem.}
\label{f:hLy001orbit}
\end{figure}

Orbits of the planar Lyapunov family are analytically computed analogously, now starting from the $E_{-1}$ equilibrium in the Hopf coordinates. These orbits get much closer to the primary than corresponding vertical ones for the same values of $L'$, and hence the effects of the perturbation are stronger and may manifest the lacks of using a lower order theory much sooner. Therefore, a wise selection of the initial conditions to propagate in the real model may be crucial to obtaining a good approximation of a planar Lyapunov orbit. As shown in Fig.~\ref{f:hLy001orbit}, if the initial conditions are taken from the analytical solution for $\ell'=0$ (left plot) the highly unstable behavior makes that the orbit very soon departs from the nominal trajectory. On the contrary, choosing $\ell=\pi$ provides a much better approximation of a periodic orbit (right plot), which, again, is easily improved by differential corrections in the original problem. 
\par

\begin{figure}[htb]
\centerline{ 
\includegraphics[scale=0.72, angle=0]{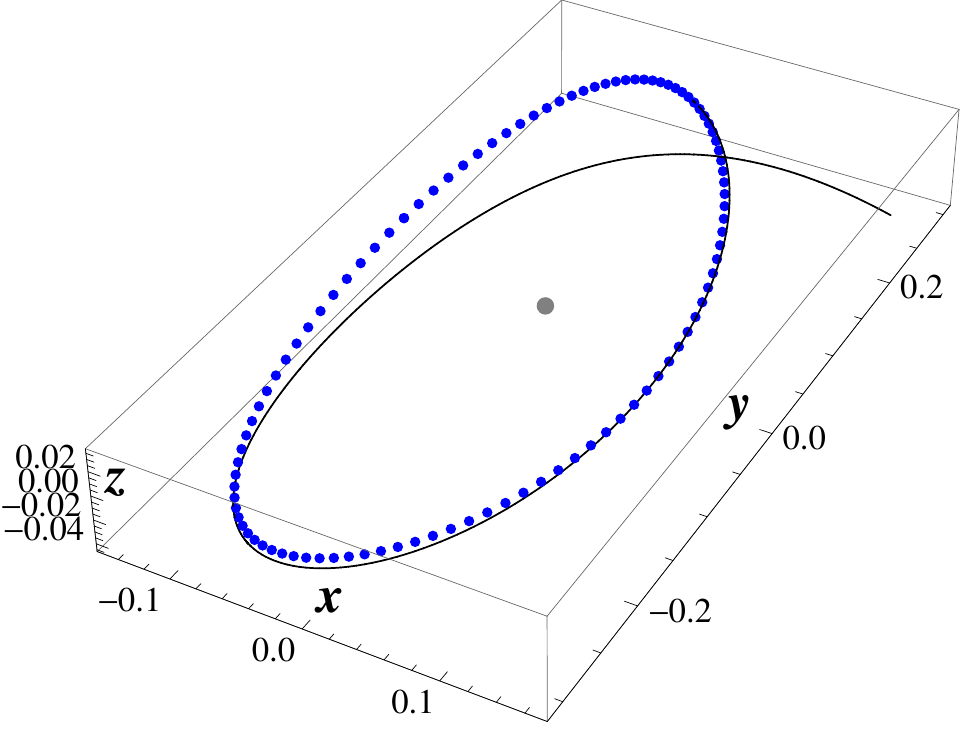}\hspace{-1cm}
\includegraphics[scale=0.7, angle=0]{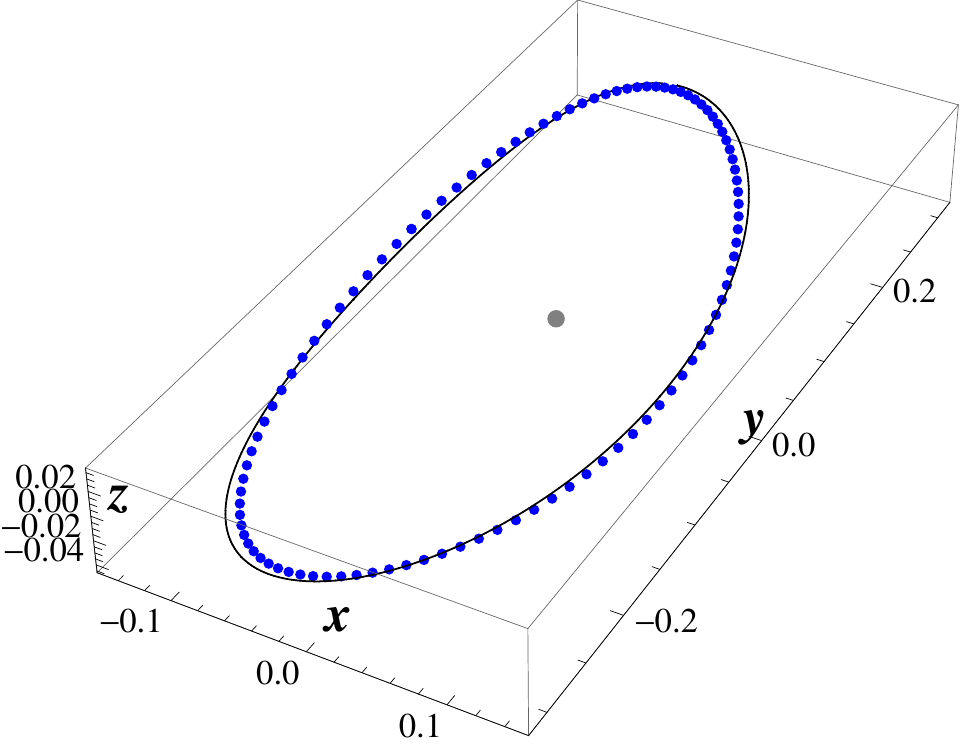}
}
\caption{Sample Halo orbit in Cartesian coordinates for $L=0.08$.}
\label{f:halorbit}
\end{figure}

The situation becomes more critic for increasing values of $L'$. Thus, while the shape of the orbit predicted by the analytical solution may resemble the aspect of the corresponding periodic orbit of the Hill problem, the initial conditions obtained from the analytical solution may fail when used as a seed in the search for a true periodic orbit. Because of the lower order of the current theory, this is exactly what happens to Halo orbits derived from the analytical solution, which the actual truncation of the theory to the second order predicts to exist only for $L'>0.077$, cf.~Eq.~(\ref{L0halo}). As shown in the left plot of Fig.~\ref{f:halorbit}, initial conditions of the analytical theory (blue points) clearly fail in providing a true Halo orbit when numerically propagated (black curve), and the situation is even worse for the hodograph (not shown). Still, the size and shape of the analytical orbit is representative of the real Halo trajectory as illustrated in the right plot of Fig.~\ref{f:halorbit}, where the black curve is a true periodic orbit of the Halo family that has been computed numerically.
\par

The situation will definitely improve when extending the analytical solution to higher orders. Nevertheless, in view of the periodic orbits of the bridge family linking vertical and planar Lyapunov orbits only exist for much higher values of $L'$, cf.~Eq.~(\ref{L12bridge}), and, therefore, can take values of the $x$ coordinate close to the Hill radius, one should not put big expectations in seizing the real dynamics of these large and highly unstable periodic orbits even when using a very high order theory.

\section{Conclusions}

A chain of five consecutive canonical transformations (a translation of the origin to the libration point, a decoupling of the linear dynamics, a polynomial normal form, the Lissajous transformation, and a short-period averaging of the elliptic anomaly) reduces the Hamiltonian of he Hill problem to a one degree of freedom Hamiltonian that comprises the dynamics in the vicinity of the libration points. The reduced phase space is advantageously described in the Hopf coordinates, which project the Lissajous variables onto the sphere and provide a deeper insight than the usual surface of section representation. In particular, the change to instability of the Lyapunov planar orbits, with the consequent bifurcation of the family of Halo orbits, is clearly visualized in the Hopf coordinates, and their associated stable and unstable manifolds are unambiguously distinguished. On the contrary, these stable and unstable manifolds almost coincide in the usual surface of section representation, in which the Lyapunov planar orbit is not represented by a point, but by the curve bounding the surface of section, and hence it is difficult to differentiate one manifold from the other.
\par

The reduced dynamics is made of fixed points of the elliptic and hyperbolic types as well as closed curves surrounding these equilibria, which are identified with the well known trajectories of the center manifold. In particular, both types of Lyapunov periodic orbits (planar and vertical) correspond to oscillations that on average remain rectilinear in Lissajous variables, whereas the Lissajous quasi-periodic orbits are, on average, rotating ellipses with variable angular momentum in Lissajous variables.
\par

The accuracy of the analytical solution is limited to energy values close to that of the libration points because of the early truncation of the developments to the fourth order in the ratio distance over Hill radius. However, the range of applicability of the solution could be trivially enlarged by the straightforward computation of higher orders of the perturbation theory.

\section*{Acknowledgemnts}

The author acknowledges partial support by the Ministry of Economic Affairs and Competitiveness of Spain, under grants ESP2013-41634-P, ESP2014-57071-R and ESP2016-76585-R.

\end{document}